\documentclass{amsart}
\usepackage{mathrsfs}
\usepackage{amsfonts}
\usepackage{amssymb}
\usepackage[papersize={6.9in,9.8in},textwidth=12.9cm,textheight=21.3cm,centering]{geometry}

\newtheorem{theorem}{Theorem}
\newtheorem{lemma}{Lemma}

\theoremstyle{definition}

\theoremstyle{remark}

\def\int{\int\limits}

\begin{document}

\title[Character sums over flat numbers]{On character sums over flat numbers }

\author{Ping Xi}

\address{School of Science, Xi'an Jiaotong University, Xi'an 710049, P. R. China} \email{xprime@163.com}

\thanks{Supported by N. S. F. (No.10601039) of P. R. China.}

\author{Yuan Yi}
\address{School of Science, Xi'an Jiaotong University, Xi'an 710049, P. R. China and Department of Mathematics, The
University of Iowa, Iowa City, Iowa 52242-1419, USA}

\email{yuanyi@mail.xjtu.edu.cn}


\begin{abstract}
Let $q\geqslant2$ be an integer, $\chi$ be any non-principal
character mod $q$, and $H=H(q)\leqslant q.$ In this paper the
authors prove some estimates for character sums of the form
\[\mathcal{W}(\chi,H;q)=\sum_{n\in\mathscr{F}(H)}\chi(n),\]where
\[\mathscr{F}(H)=\left\{n\in\mathbb{Z}\bigg|(n,q)=1,1\leqslant n,\overline{n}\leqslant q,
|n-\overline{n}|\leqslant H\right\},\] $\overline{n}$ is defined by
$n\overline{n}\equiv 1\pmod q.$
\end{abstract}

 \maketitle

% \tableofcontents

\vspace*{-1cm}

\section{Introduction}

Let $q\geqslant2$ be an integer, $\chi$ be a non-principal character
mod $q$. It is quite an important problem in analytic number theory
to obtain a sharp estimate for the character sum
\[\sum_{x=N+1}^{N+H}\chi(f(x)),\] where $f(x)\in\mathbb{Z}[x]$, $N$
and $H$ are positive integers. The classical result, due to
P\'{o}lya and Vinogradov \cite{P,V}, is the
estimate\[\sum_{n=N+1}^{N+H}\chi(n)\ll q^{1/2}\log q,\]where $\ll$
is the Vinogradov's notation. About half a century later, Burgess'
immortal work \cite{B1,B2} showed that
\[\sum_{n=N+1}^{N+H}\chi(n)\ll
H^{1-1/r}q^{(r+1)/4r^2+o(1)}\]holds with $r=1,2,3$ for any $q$ and
with arbitrary positive integer $r$ if $q$ is cube-free. Under
Generalized Riemann Hypothesis, Montgomery and Vaughan \cite{MV}
sharpened the P\'{o}lya-Vinogradov bound to
\[\sum_{n=N+1}^{N+H}\chi(n)\ll q^{1/2}\log\log q.\]

Burgess' estimate is such a milestone that nobody can
unconditionally beat the barrier in general by any advanced
technology (Some partial improvements and progresses can be found in
\cite{GS1,GS2,G}, et al). However, estimates for the character sums
over special numbers and sequences, such as factorials, Beatty
sequences, binomial coefficients and other combinatorial numbers,
have attracted many scholars' interests. A complete list of the
results and open problems are referred to \cite{S}.

In this paper, we shall deal with another kind of special numbers.
It is known that the distribution of $\overline{n}$ is quite
irregular, where $\overline{n}$ is the inverse of $n$ mod $q$, i.e.
$n\overline{n}\equiv1\pmod q.$ How about the distribution of
$|n-\overline{n}|$? In \cite{Z}, W. Zhang proved that
\begin{equation}\label{eq:1}\sideset{}{^*}\sum_{\substack{n=1\\|n-\overline{n}|\leqslant\delta q}}^q1=\delta(2-\delta)\varphi(q)+O\left(q^{1/2}\tau^2(q)\log^3q\right),
\end{equation}
where $\delta\in(0,1]$ is a constant, $\varphi(q)$ is the Euler
function and $\tau(q)$ is the divisor function, $\sum^*$ denotes the
summation over the integers that are coprime to $q$.

In fact, W. Zhang studied the number of the integers that are within
a given distance to their inverses mod $q$. Now we consider the
character sums over these integers. We shall present the problem of
a more general case.

Let $q\geqslant2$ be a fixed integer and $H=H(q)\leqslant q.$ We put
\[\mathscr{F}(H)=\left\{n\in\mathbb{Z}\bigg|(n,q)=1,1\leqslant n,\overline{n}\leqslant q,
|n-\overline{n}|\leqslant H\right\}.\] Each element in
$\mathscr{F}(H)$ is called a $H$-flat number mod $q$. Note that in
the definition of $\mathscr{F}(H)$, the size of $H$ is $O(q),$ not
necessary being $H\asymp q$ as in (\ref{eq:1}).

In this paper, we shall study the character sums over such $H$-flat
numbers mod $q$. That is we shall prove nontrivial upper bounds for
\begin{equation}\label{eq:2}\mathcal{W}(\chi,H;q)=\sum_{n\in\mathscr{F}(H)}\chi(n).\end{equation}
It is obvious that $n\in\mathscr{F}(H)$ implies
$q-n\in\mathscr{F}(H)$, thus $\chi(n)+\chi(q-n)=0$ if $\chi(-1)=-1,$
so $\mathcal{W}(\chi,H;q)=0$. Hence we only deal with the case with
$\chi(-1)=1$ throughout this paper.
\begin{theorem}\label{thm:1}Let $q\geqslant2$, $\chi$ be a non-principal
character mod $q.$ Then  we have
\[\mathcal{W}(\chi,H;q)\ll q^{1/2}\tau^2(q)\log H.\]\end{theorem}

The proof of Theorem \ref{thm:1} depends on the estimate for the
general Kloosterman sums twisted by Dirichlet characters, and the
upper bound in Theorem \ref{thm:1} is independent of $H$, to be
precise, the result may be trivial if $H$ is quite small. However if
$q$ is odd, and $\chi$ is the Jacobi symbol mod $q$ (which reduces
to Legendre symbol if $q$ is a prime), we have corresponding
calculation formulae for this Kloosterman sum, known as Sali\'{e}
sum, and we can obtain an upper bound depending on $H$, which can be
stated as follows
\begin{theorem}\label{thm:2}Let $q\geqslant3$ be an odd square-free integer, $\chi$ be the Jacobi symbol mod $q.$ Then we have
\[\mathcal{W}(\chi,H;q)\ll H^{1-1/r}q^{(r+1)/4r^2}\tau(q)\log q,\]
where $r\geqslant1$ is an arbitrary integer. \end{theorem}

\bigskip

\section{General Kloostermann sums and character sums}
The classical Kloosterman sum is defined by
\[S(m,n;q)=\sideset{}{^*}\sum_{a\bmod
q}e\left(\frac{ma+n\overline{a}}{q}\right),\]where $e(x)=e^{2\pi
ix}$. The well-known upper bound essentially due to A. Weil \cite{W}
is
\[S(m,n;q)\ll q^{1/2}(m, n,
q)^{1/2}\tau(q),\]where $(m,n,q)$ denotes the greatest common
divisor of $m,n,q$.

In the proof of following sections, we require a general Kloosterman
sum twisted by a Dirichlet character such as
\[S_\chi(m,n;q)=\sideset{}{^*}\sum_{a\bmod
q}\chi(a)e\left(\frac{ma+n\overline{a}}{q}\right).\] Taking
$\chi=\chi^0$ as the principal character mod $q$, this reduces to
$S(m,n;q)$.

We require an upper bound estimation for $S_\chi(m,n;q)$, the
original proofs \cite{W,E} carry over with minor modifications.
\begin{lemma}\label{lm:1}
Let $q$ be a positive integer, then we
have\[S_\chi(m,n;q)\ll q^{1/2}(m, n,
q)^{1/2}\tau(q).\]
\end{lemma}

\begin{lemma}\label{lm:2} Let $q\geqslant2$, $\chi$ be a Dirichlet character mod $q$. For any $d$ with $d|q$, we
define
\[T_\chi(m,n;d,q)=\sideset{}{^*}\sum_{a\bmod
q}\chi(a)e\left(\frac{ma+n\overline{a}}{d}\right),\]where
$a\overline{a}\equiv1(\bmod q)$. Then for any $d$ with
$d\ell=q,(d,\ell)=1$, we have
\[T_\chi(m,n;d,q)=
\begin{cases}\varphi(\ell)S_{\chi_1}(m,n;d),& \text{if\ \ }\chi_2=\chi_2^0,\\
0,& \text{if\ \ }\chi_2\neq\chi_2^0,
\end{cases}\]
where $\chi_1\bmod d,\chi_2\bmod \ell$ with $\chi_1\chi_2=\chi,$ and
$\chi_2^0$ is the principal character mod $\ell$.
\end{lemma}

\proof Let $a=a_1\ell+a_2d,$ then
\[\begin{split}T_\chi(m,n;d,q)&=\sideset{}{^*}\sum_{a_1\bmod
d}~~\sideset{}{^*}\sum_{a_2\bmod
\ell}\chi(a_1\ell+a_2d)e\left(\frac{ma_1\ell+n\overline{a_1\ell}}{d}\right)\\
&=\sideset{}{^*}\sum_{a_1\bmod d}~~\sideset{}{^*}\sum_{a_2\bmod
\ell}\chi_1(a_1\ell)\chi_2(a_2d)e\left(\frac{ma_1\ell+n\overline{a_1\ell}}{d}\right)\\
&=\chi_2(d)\sideset{}{^*}\sum_{a_1\bmod
d}\chi_1(a_1\ell)e\left(\frac{ma_1\ell+n\overline{a_1\ell}}{d}\right)\sideset{}{^*}\sum_{a_2\bmod
\ell}\chi_2(a_2),\end{split}\] then the lemma follows from the
orthogonality of Dirichlet characters.  \endproof

If $q$ is odd, and $\chi$ is the Jacobi symbol mod $q$, we have a
calculation formula of the general Kloosterman sums, known as
Sali\'{e} sums (See \cite{I}, Lemma 4.9).
\begin{lemma}\label{lm:3} If $(q,2n)=1,$ and $\chi$ is the Jacobi symbol mod $q$, then we have
\[S_\chi(m,n;q)=
\varepsilon_qq^{1/2}\chi(n)\sum_{y^2\equiv mn(\bmod
q)}e\left(\frac{2y}{q}\right),\] where $\varepsilon_q$ is a constant
with $|\varepsilon_q|=1$.
\end{lemma}

We also require Burgess' classical result on character sums, see
\cite{B1}, Theorem 2.
\begin{lemma}\label{lm:4}
If $q\geqslant2$ is square-free, then for any non-principal
character $\chi$ mod $q$, we have\[\sum_{n=N+1}^{N+A}\chi(n)\ll
A^{1-1/r}q^{(r+1)/4r^2}\log q,\] where $r\geqslant1$ is an arbitrary
integer.
\end{lemma}

\bigskip

\section{Proof of Theorem \ref{thm:1}}
It obvious that $\mathcal{W}(\chi,H;q)$ has the same essential bound
with
\[\mathcal{W}^*(\chi,H;q)=\sum_{t\leqslant H}\sideset{}{^*}
\sum_{\substack{n\leqslant q\\n-\overline{n}\equiv t(\bmod
q)}}\chi(n).\]

We denote $\|x\|=\min_{n\in\mathbb{Z}}|x-n|$. Apply the
identity\[\sum_{n=1}^qe\left(\frac{an}{q}\right)=
\begin{cases}
q,& q\mid a,\\
0,& q\nmid a,
\end{cases}\]
together with Lemma \ref{lm:1} we can obtain that
\[\begin{split}\mathcal{W}^*(\chi,H;q)&=\frac{1}{q}\sum_{m\leqslant q}\sum_{t\leqslant H}e\left(-\frac{mt}{q}\right)
\sideset{}{^*}\sum_{n\leqslant q}\chi(n)e\left(m\frac{n-\overline{n}}{q}\right)\\
&\ll\frac{1}{q}\sum_{m\leqslant
q-1}\min\Big(H,\Big\|\frac{m}{q}\Big\|^{-1}\Big)\left|
S_\chi(m,-m;q)\right|\\
 &\ll
q^{-1/2}\tau(q)\sum_{m\leqslant q-1}(m,q)^{1/2}\min\Big(H,\Big\|\frac{m}{q}\Big\|^{-1}\Big)\\
&\ll
q^{-1/2}\tau(q)\sum_{m\leqslant q-1}(m,q)^{1/2}\min\Big(H,\frac{q}{m}\Big)\\
&\ll Hq^{-1/2}\tau(q)\sum_{m\leqslant q/H}(m,q)^{1/2}+
q^{1/2}\tau(q)\sum_{q/H<m\leqslant
q-1}\frac{(m,q)^{1/2}}{m}.\end{split}\]

By the following calculations,
\[\sum_{m\leqslant q/H}(m,q)^{1/2}\ll\sum_{d|q}d^{1/2}\sum_{\substack{m\leqslant q/H\\ d|m}}1\ll H^{-1}q\tau(q),\] and
\[\begin{split}\sum_{q/H<m\leqslant
q-1}\frac{(m,q)^{1/2}}{m}&\ll\sum_{d|q}d^{1/2}\sum_{\substack{q/H<m\leqslant
q-1\\ d|m}}\frac{1}{m}\\
&\ll \sum_{d|q}d^{-1/2}\sum_{q/Hd<m\leqslant q/d}\frac{1}{m}\\
&\ll \tau(q)\log H,\end{split}\] we have
\[\mathcal{W}^*(\chi,H;q)\ll q^{1/2}\tau^2(q)\log H.\] And $\mathcal{W}(\chi,H;q)$ has the
same bound. This completes the proof of Theorem \ref{thm:1}.

\bigskip

\section{Proof of Theorem \ref{thm:2}}

In this section, we shall deal with a special case of (\ref{eq:2}), that is
$q$ being an odd square-free integer, and $\chi=(\frac{\cdot}{q})$
being the Jacobi symbol mod $q$.

Following the similar arguments in Section 3, we have
\[\begin{split}\mathcal{W}^*(\chi,H;q)&=\frac{1}{q}\sum_{m\leqslant q}\sum_{t\leqslant H}e\left(-\frac{mt}{q}\right)
S_\chi(m,-m;q)\\
&=\frac{1}{q}\sum_{d|q}\sum_{\substack{m=1\\(m,q)=d}}^{q}\sum_{t\leqslant
H}e\left(-\frac{mt}{q}\right) S_\chi(m,-m;q)\\
&=\frac{1}{q}\sum_{d|q}\sum_{\substack{m=1\\(m,q/d)=1}}^{q/d}\sum_{t\leqslant
H}e\left(-\frac{mdt}{q}\right) S_\chi(md,-md;q)\\
&=\frac{1}{q}\sum_{d|q}\sideset{}{^*}\sum_{m\leqslant
d}\sum_{t\leqslant H}e\left(-\frac{mt}{d}\right)
\sideset{}{^*}\sum_{a\bmod
q}\chi(a)e\left(\frac{ma-m\overline{a}}{d}\right)\\
&=\frac{1}{q}\sum_{d|q}\sideset{}{^*}\sum_{m\leqslant
d}\sum_{t\leqslant H}e\left(-\frac{mt}{d}\right)
T_\chi(m,-m;d,q).\end{split}\]

We write $q=d\ell,$ where $(d,\ell)=1$ since $q$ is square-free. We
also write $\chi_1\bmod d,\chi_2\bmod \ell$ with
$\chi_1\chi_2=\chi.$ Note that $\chi$ is a real primitive character
mod $q$, so from Lemma \ref{lm:2} we know that, $\chi_1$ and
$\chi_2$ must be real primitive characters mod $d$ and $\ell$
respectively, thus $T_\chi(m,-m;d,q)=0$ if $\ell>1$. Applying Lemma
\ref{lm:2} and Lemma \ref{lm:3}, we can deduce that
\[\begin{split}\mathcal{W}^*(\chi,H;q)&=\frac{1}{q} \sideset{}{^*}\sum_{m\leqslant q}\sum_{t\leqslant
H}e\left(-\frac{mt}{q}\right)S_\chi(m,-m;q)\\
&=\frac{\varepsilon_q}{q^{1/2}}\sum_{t\leqslant H}\sum_{m\leqslant
q} \chi(m)e\left(-\frac{mt}{q}\right)\sum_{y^2\equiv -m^2(\bmod
q)}e\left(\frac{2y}{q}\right)\\
&=\frac{\varepsilon_q}{q^{1/2}}\sum_{t\leqslant H}\sum_{m\leqslant
q} \chi(m)e\left(-\frac{mt}{q}\right)\sum_{\delta^2\equiv -1(\bmod
q)}e\left(\frac{2\delta m}{q}\right).\end{split}\] Thus
\[\begin{split}\mathcal{W}^*(\chi,H;q)
&=\frac{\varepsilon_q}{q^{1/2}}\tau(\chi)\sum_{\delta^2\equiv
-1(\bmod q)}\sum_{t\leqslant H}\chi(2\delta-t),\end{split}\] where
$\tau(\chi)=\sum_{n\bmod q}\chi(n)e(\frac{n}{q})$ is the Gauss sum.

Applying Lemma \ref{lm:4} and $|\tau(\chi)|=q^{1/2}$, we obtain that
\[\mathcal{W}^*(\chi,H;q)\ll H^{1-1/r}q^{(r+1)/4r^2}\tau(q)\log q\] for any integer $r\geqslant1$.
And $\mathcal{W}(\chi,H;q)$ has the same bound. This completes the
proof of Theorem \ref{thm:2}.

\bigskip

\section{Further discussion}

Each positive integer can be represented as the product of two
coprime parts, one is square-free and the other is square-full. In
fact, the method in Section 4 can lead to a nontrivial estimate for the
modulo $q$, whose square-full part is quite small in comparison with
$q$. This estimate mainly depends on the following lemma. Henceforth
the undefined symbols and notations are all referred to this lemma.
\begin{lemma}\label{lm:5} Let $q\geqslant2$, $\chi$ be a Dirichlet character mod $q$.
We write $q=q_1q_2,$ where $q_1$ is the square-free part, and $q_2$
is square-full. For any $d$ with $d|q$, we write $d=d_1d_2,$ where
$d_1|q_1$ and $d_2|q_2$. Then
\[T_\chi(m,n;d,q)=
\begin{cases}
\chi_1(q_2)\chi_2(q_1)\varphi(\ell_2)S_{\chi'}(m\ell_2,n\overline{q_2^2}\ell_2;d_1)\\
\ \ \ \ \times T_{\chi_2}(m\ell_1,n\overline{q_1^2}\ell_1;d_2,q_2),& \text{if\ \ }\chi''=\chi''_0,\\
0,& \text{if\ \ }\chi''\neq\chi''_0,
\end{cases}\]
where $q_1\overline{q_1}\equiv1(\bmod
q_2),q_2\overline{q_2}\equiv1(\bmod q_1),$ $\chi_1\bmod
q_1,\chi_2\bmod q_2$ with $\chi_1\chi_2=\chi,$ and
 $\chi'\bmod d_1,\chi''\bmod q_1/d_1$ with $\chi'\chi{''}=\chi_1,$ $\chi''_0$ is the principal character mod $q_1/d_1.$
\end{lemma}

\proof For convenience, we write $\ell_1=q_1/d_1,\ell_2=q_2/d_2.$
Let $a=a_1q_2+a_2q_1,$ then
\begin{equation}\label{eq:3}\begin{split}T_\chi(m,n;d,q)
&=\sideset{}{^*}\sum_{a_1\bmod q_1}~~\sideset{}{^*}\sum_{a_2\bmod
q_2}\chi(a_1q_2+a_2q_1)e\left(\frac{m(a_1q_2+a_2q_1)+n\overline{(a_1q_2+a_2q_1)}}{d_1d_2}\right)\\
&=\sideset{}{^*}\sum_{a_1\bmod q_1}~~\sideset{}{^*}\sum_{a_2\bmod
q_2}\chi(a_1q_2+a_2q_1)e\left(\frac{m(a_1q_2+a_2q_1)+n(\overline{a_1q_2^2}q_2+\overline{a_2q_1^2}q_1)}{d_1d_2}\right)\\
&=\chi_1(q_2)\chi_2(q_1)\sideset{}{^*}\sum_{a_1\bmod q_1}\chi_1(a_1)e\left(\frac{m\ell_2a_1+n\overline{q_2^2}\ell_2\overline{a_1}}{d_1}\right)\\
&\ \ \ \ \times\sideset{}{^*}\sum_{a_2\bmod
q_2}\chi_2(a_2)e\left(\frac{m\ell_1a_2+n\overline{q_1^2}\ell_1\overline{a_2}}{d_2}\right)\\
&=\chi_1(q_2)\chi_2(q_1)T_{\chi_1}(m\ell_2,n\overline{q_2^2}\ell_2;d_1,q_1)T_{\chi_2}(m\ell_1,n\overline{q_1^2}\ell_1;d_2,q_2).\end{split}\end{equation}

Now we write $\chi_1=\chi'\chi{''},$ where $\chi'\bmod
d_1,\chi''\bmod \ell_1.$ From Lemma \ref{lm:4} and (\ref{eq:3}), we
can deduce this lemma.
\endproof

Now we would like to give a slight improvement of Theorem
\ref{thm:2}.

First, by Lemma \ref{lm:5} and the method in Section 4, we have
\[\begin{split}\mathcal{W}^*(\chi,H;q)&=\frac{1}{q}\sum_{d_1|q_1}\sum_{d_2|q_2}\sideset{}{^*}\sum_{m\leqslant d_1d_2}\sum_{t\leqslant
H}e\left(-\frac{mt}{d_1d_2}\right) T_\chi(m,-m;d_1d_2,q)\\
&=\frac{\chi_1(q_2)\chi_2(q_1)}{q}\sum_{\substack{d_1|q_1\\
\chi''=\chi''_0}}\sum_{d_2|q_2}\varphi(\ell_2)\sideset{}{^*}\sum_{m\leqslant
d_1d_2}\sum_{t\leqslant H}e\left(-\frac{mt}{d_1d_2}\right)
\\
&\ \ \ \ \times
S_{\chi'}(m\ell_2,-m\overline{q_2^2}\ell_2;d_1)T_{\chi_2}(m\ell_1,-m\overline{q_1^2}\ell_1;d_2,q_2).\end{split}\]

Now opening the Kloosterman sums, we obtain
\begin{eqnarray}\mathcal{W}^*(\chi,H;q)&=&\frac{\chi_1(q_2)\chi_2(q_1)}{q}\sum_{\substack{d_1|q_1\\
\chi''=\chi''_0}}\sum_{d_2|q_2}\varphi(\ell_2)\sum_{t\leqslant
H}\sideset{}{^*}\sum_{a\leqslant q_2}\chi_2(a)\nonumber
\\
& &\ \ \ \ \times\sideset{}{^*}\sum_{m\leqslant
d_1d_2}e\left(-\frac{mt}{d_1d_2}\right)
e\left(\frac{m\ell_1a-m\overline{q_1^2}\ell_1\overline{a}}{d_2}\right)S_{\chi'}(m\ell_2,-m\overline{q_2^2}\ell_2;d_1)\nonumber \\
&\ll& q^{-1}\sum_{\substack{d_1|q_1\\
\chi''=\chi''_0}}d_1^{1/2}\sum_{d_2|q_2}\varphi(\ell_2)\sideset{}{^*}\sum_{a\leqslant
q_2}\Bigg|\sum_{t\leqslant H}~~\sideset{}{^*}\sum_{m\leqslant
d_1d_2}\chi'(m)\nonumber
\\
& &\ \ \ \ \times
e\left(\frac{mt}{d_1d_2}\right)e\left(\frac{m\overline{q_1^2}\ell_1\overline{a}-m\ell_1a}{d_2}\right)
\sum_{y^2\equiv -m^2\overline{q_2^2}\ell_2^2(\text{mod}
d_1)}e\left(\frac{2y}{d_1}\right)\Bigg|\nonumber \\
&=&q^{-1}\sum_{\substack{d_1|q_1\\
\chi''=\chi''_0}}d_1^{1/2}\sum_{d_2|q_2}\varphi(\ell_2)\sideset{}{^*}\sum_{a\leqslant
q_2}\Bigg|\sum_{t\leqslant H}~~\sideset{}{^*}\sum_{m\leqslant
d_1d_2}\chi'(m)\nonumber
\\
& &\ \ \ \ \times
e\left(\frac{mt}{d_1d_2}\right)e\left(\frac{m\overline{q_1^2}\ell_1\overline{a}-m\ell_1a}{d_2}\right)
\sum_{\eta^2\equiv-1(\bmod d_1)}e\left(\frac{2\eta
m\overline{q_2}\ell_2}{d_1}\right)\Bigg|\nonumber \\
&=&q^{-1}\sum_{\substack{d_1|q_1\\
\chi''=\chi''_0}}d_1^{1/2}\sum_{d_2|q_2}\varphi(\ell_2)\sideset{}{^*}\sum_{a\leqslant
q_2}\Bigg|\sum_{\eta^2\equiv-1(\bmod d_1)}\sum_{t\leqslant
H}\nonumber
\\
& &\ \ \ \ \times\sideset{}{^*}\sum_{m\leqslant
d_1d_2}\chi'(m)e\left(m\frac{t+\overline{q_1^2}\ell_1d_1\overline{a}-\ell_1d_1a+2\eta
\overline{q_2}\ell_2d_2}{d_1d_2}\right)\Bigg|\nonumber \\
&\ll& q^{-1}\sum_{\substack{d_1|q_1\\
\chi''=\chi''_0}}d_1^{1/2}\sum_{\eta^2\equiv-1(\bmod
d_1)}\sum_{d_2|q_2}\varphi(\ell_2)\sideset{}{^*}\sum_{a\leqslant q_2}\nonumber \\
& &\ \ \ \ \times\Bigg|\sum_{t\leqslant
H}~~\sideset{}{^*}\sum_{m\leqslant
d_1d_2}\chi'\widetilde{\chi}(m)e\left(m\frac{t+\overline{a}-a+2\eta}{d_1d_2}\right)\Bigg|\nonumber \\
&=&q^{-1}\sum_{\substack{d_1|q_1\\
\chi''=\chi''_0}}d_1^{1/2}\sum_{\eta^2\equiv-1(\text{mod}
d_1)}\sum_{d_2|q_2}\varphi(\ell_2)\sideset{}{^*}\sum_{a\leqslant
q_2}\left|\sum_{t\leqslant
H}G(t+\overline{a}-a+2\eta,\chi'\widetilde{\chi})\right|\label{eq:4},\end{eqnarray}
here $\widetilde{\chi}$ is the principal character mod $d_2,$ and
\[G(s,\chi)=\sideset{}{^*}\sum_{m=1}^{q}\chi(m)e\left(\frac{ms}{q}\right)\]is
the Gauss sum corresponding to $\chi$ mod $q$.

Note that $\chi_1\chi_d^0$ is not a primitive character mod $dq_1,$
thus $G(s,\chi)$ cannot be expressed in terms of $\tau(\chi)$
directly. Thus we would like to refer to the following lemma.
\begin{lemma}\label{lm:6} Let $q\geqslant2$, $\chi$ be a non-principal character mod
$q$. We write $q=q_1q_2,$ $(q_1,q_2)=1,$ and $\chi_1\bmod
q_1,\chi_2\bmod q_2$ with $\chi_1\chi_2=\chi.$

If $q_1$ is square-free, $\chi_1$ is a primitive character mod $q_1$
and $\chi_2$ is the principal character mod $q_2$, then arbitrary
integer $r\geqslant1$, we have
\[\sum_{n=M+1}^{M+N}G(n,\chi)\ll(Nq_2^{-1})^{1-1/r}q_1^{1/2+(r+1)/4r^2}\varphi(q_2)2^{\omega(q_2)}\log
q_1,\] where $\omega(q_2)$ denotes the number of distinct prime
factors of $q_2$.
\end{lemma}

\proof First, Gauss sums enjoy the following properties(see
\cite{PP}, Lemma 1.1, 1.2):
\[G(n,\chi)=\chi_1(q_2)\chi_2(q_1)G(n,\chi_1)G(n,\chi_2),\]
and
\[G(n,\chi_2)=\sideset{}{^*}\sum_{m=1}^{q_2}e\left(\frac{mn}{q_2}\right)=\mu\left(\frac{q_2}{(n,q_2)}\right)\varphi(q_2)
\varphi^{-1}\left(\frac{q_2}{(n,q_2)}\right),\] which give
\[G(n,\chi)=\chi_1(\overline{n}q_2)\tau(\chi_1)\mu\left(\frac{q_2}{(n,q_2)}\right)\varphi(q_2)
\varphi^{-1}\left(\frac{q_2}{(n,q_2)}\right).\]

Thus we can deduce that
\[\begin{split}\sum_{n=M+1}^{M+N}G(n,\chi)&=\chi_1(q_2)\tau(\chi_1)\varphi(q_2)\sum_{n=M+1}^{M+N}\overline{\chi_1}(n)\mu\left(\frac{q_2}{(n,q_2)}\right)
\varphi^{-1}\left(\frac{q_2}{(n,q_2)}\right)\\
&=\chi_1(q_2)\tau(\chi_1)\varphi(q_2)\sum_{d|q_2}\sum_{\substack{n=M+1\\(n,q_2)=d}}^{M+N}\overline{\chi_1}(n)\mu\left(\frac{q_2}{d}\right)
\varphi^{-1}\left(\frac{q_2}{d}\right)\\
&=\chi_1(q_2)\tau(\chi_1)\varphi(q_2)\sum_{d|q_2}\mu\left(\frac{q_2}{d}\right)
\varphi^{-1}\left(\frac{q_2}{d}\right)\sum_{\substack{n=M+1\\(n,q_2)=d}}^{M+N}\overline{\chi_1}(n)\\
&\ll
q_1^{1/2}\varphi(q_2)\sum_{d|q_2}\mu^2\left(\frac{q_2}{d}\right)
\varphi^{-1}\left(\frac{q_2}{d}\right)\left|\sum_{\substack{M/d<n\leq
(M+N)/d\\(n,q_2/d)=1}}\overline{\chi_1}(n)\right|\\
&\ll
q_1^{1/2}\varphi(q_2)\sum_{d|q_2}\mu^2\left(\frac{q_2}{d}\right)\tau\left(\frac{q_2}{d}\right)
\varphi^{-1}\left(\frac{q_2}{d}\right)\left(\frac{N}{d}\right)^{1-1/r}q_1^{(r+1)/4r^2}\log
q_1\\
&\ll (Nq_2^{-1})^{1-1/r}q_1^{1/2+(r+1)/4r^2}\varphi(q_2)\log
q_1\sum_{d|q_2}\mu^2(d)\tau(d)\varphi^{-1}(d)d^{1-1/r}\\
&\ll
(Nq_2^{-1})^{1-1/r}q_1^{1/2+(r+1)/4r^2}\varphi(q_2)2^{\omega(q_2)}\log
q_1,\end{split}\] where $\omega(q_2)$ denotes the number of distinct
prime factors of $q_2$ and $r\geqslant1$ is an arbitrary integer.
This completes the proof of Lemma \ref{lm:6}.
\endproof

From (\ref{eq:4}) and Lemma \ref{lm:6}, we obtain
\[\begin{split}\mathcal{W}^*(\chi,H;q)
&\ll
q^{-1}\varphi(q_2)H^{1-1/r}\sum_{d_1|q_1}d_1^{1+(r+1)/4r^2}\tau(d_1)\log
d_1\\
&\ \ \ \ \times\sum_{d_2\ell_2=q_2}\varphi(\ell_2)\varphi(d_2)d_2^{1/r}2^{\omega(d_2)}\\
&\ll(Hq_2^{-1})^{1-1/r}\varphi^2(q_2)q_1^{(r+1)/4r^2}2^{\omega(q)}\log
q_1.\end{split}\]

In summary, the discussions above can be stated as the final
theorem.
\begin{theorem}\label{thm:3}Let $q\geqslant3$ be an odd integer, $\chi$ be the Jacobi symbol mod $q.$ Then we have
\[\mathcal{W}(\chi,H;q)\ll (Hq_2^{-1})^{1-1/r}\varphi^2(q_2)q_1^{(r+1)/4r^2}2^{\omega(q)}\log
q_1,\] where $q_1$ is the square-free part of $q$, $q_2$ is
square-full, $\omega(q)$ denotes the number of distinct prime
factors of $q$, and $r\geqslant1$ is an arbitrary integer.
\end{theorem}

We should point out again that Theorem \ref{thm:3} is nontrivial when the
square-full part of $q$ is quite small in comparison with $q$.

\bigskip

\bigskip

\bibliographystyle{plainnat}

\end{document}